\documentclass[11pt]{amsart}
\usepackage{mathptmx}
\usepackage{amsmath}
\usepackage{amscd}
\usepackage{amsfonts}
\usepackage{amssymb}
\usepackage{graphicx,cite}
\usepackage{mathrsfs}
\usepackage[all,cmtip,2cell]{xy}

\newtheorem{dfn}{Definition}[section]
\newtheorem{thm}{Theorem}
\newtheorem{lem}[thm]{Lemma}
\newtheorem {pro}[thm]{Proposition}

\title{Cross-connections of linear transformation semigroup}
\author{P. A. Azeef Muhammed}
\address{School of Mathematics, Indian Institute of Science Education and Research, Thiruvananthapuram, Kerala-695016, India.\footnote{Present address: Institute of Mathematics and Computer Science, Ural Federal University, Lenina 51, 620000 Ekaterinburg, Russia. E-mail address: a.a.parail@urfu.ru}}
\email{azeefp@gmail.com}
\thanks{The author acknowledges the financial support of IISER, Thiruvananthapuram in the preparation of this article.
On the final stages, his work was also supported by the Competitiveness Enhancement Program of Ural Federal University, Russia.}  
\date{}

\begin{document}

\begin{abstract}
Cross-connection theory developed by Nambooripad  is the construction of a semigroup from its principal left (right) ideals using categories. We briefly describe the general cross-connection theory for regular semigroups and use it to study the {normal categories} arising from the semigroup $Sing(V)$ of singular linear transformations on an arbitrary vectorspace $V$ over a field $K$. There is an inbuilt notion of duality in the cross-connection theory, and we observe that it coincides with the conventional algebraic duality of vector spaces. We describe various cross-connections between these categories and show that although there are many cross-connections, upto isomorphism, we have only one semigroup arising from these categories. But if we restrict the categories suitably, we can construct some interesting subsemigroups of the {variant} of the linear transformation semigroup. 

\end{abstract}

\maketitle

\section{Background and overview}
We begin with a brief informal discussion on the general theory of cross-connections which will help us to place our results in a proper context. In the study of the structure theory of regular semigroups, Hall \cite{hall} used the ideals of a given regular semigroup to analyse its structure. In 1974, Grillet \cite{gril,gril1,gril2} refined Hall's theory to abstractly characterize the sets of the principal ideals of a given regular semigroup as \emph{regular} partially ordered sets, and constructed the fundamental image of the given semigroup as a {cross-connection} semigroup.

Grillet explicitly described the relationship between the principal ideals of a given fundamental regular semigroup and showed that they are {cross-connected} with each other via two order preserving mappings $\Gamma$ and $\Delta$, which he rightly called {cross-connections}. He also showed that given two {appropriate} regular partially ordered sets and a pair of cross-connections, the cross-connection semigroup obtained is a fundamental regular semigroup. It can be seen that given a regular semigroup and its {full} regular subsemigroup, both induces the same cross-connection (upto isomorphisms), and hence we get the same fundamental semigroup. In other words, isomorphic {biordered sets} give rise to isomorphic cross-connections \cite{bicxn}. Thus by Grillet's construction using partially ordered sets, we are only able to capture the {skeleton} of the semigroup, not the complete semigroup.

Nambooripad recognized this deficiency and remarked that the data provided by partially ordered sets are insufficient to completely characterize arbitrary regular semigroups. He observed that any partially ordered set can be viewed as a {preorder category}, and using the richer structure provided by {small} categories (treating them as algebraic objects), one can in fact generalize Grillet's theory. In 1994, Nambooripad \cite{cross}\nocite{mem} constructed arbitrary regular semigroups as cross-connection semigroups by characterizing the principal ideals of a regular semigroup as \emph{normal categories} and replacing the cross-connection maps by a pair of functors. Although philosophically the construction was very much same as that of Grillet's, Nambooripad's work was a big leap in terms of technique, that Meakin and Rajan \cite{jmar} in their tribute to Nambooripad described the theory of cross-connections as a \emph{technical tour de force}.

A {normal category} is a categorical abstraction of the partially ordered set of principal left (right) ideals of a regular semigroup $S$. To replace partially ordered sets of Grillet with categories, morphisms were to be introduced and the natural choice for morphisms were mappings induced by partial translations. Hence a normal category is a small category whose object set is a regular partially ordered set with morphisms as mappings between them. 

The morphisms of a normal category has an important factorization property called \emph{normal factorization} by which each translation could be factored into three components---a {retraction}, an {isomorphism} and an {inclusion}. This factorization results from the ideal structure of the regular semigroup. Given a translation $f$ between the principal ideals $I$ and $J$, by the {retraction} $q$, one can \emph{go down} to an appropriate ideal $I'\subseteq I$; so that the {isomorphism} $u$, will translate it to an ideal $J'\subseteq J$ where $J'= \text{Im} f$. Then the {inclusion} $j = j(J',J)$, is the inclusion map from $J'$ to $J$ so that $f=quj$. Whenever $I'\subseteq I$, we can see that there is a retraction $q: I\to I'$ such that $j(I',I)q =1_{I'}$. This property of the inclusion morphisms is called \emph{splitting}.  

Observe that the inclusion morphisms are the inclusions of the regular partially ordered set, and hence the partially ordered set $P$ may be identified as a strict-preorder subcategory $\mathcal{P}$ (whose morphisms are precisely the inclusions) of the normal category. Thus we have a category $\mathcal{C}$ with a \emph{choice of subobjects} where the regular partially ordered set is sitting inside as the subcategory $\mathcal{P}$, and hence denoted as $(\mathcal{C},\mathcal{P})$.

A 'nice' cluster of morphisms with a common codomain is referred to as a {normal cone}. This is a generalization of the {normal mapping} of Grillet and they act as the basic building blocks of our construction. Given a partially ordered set $P$, a normal mapping with {apex} $a \in P$, was defined by Grillet as an order preserving mapping from $P$ to the principal order ideal $P(a)$. Given a category $(\mathcal{C},\mathcal{P})$ with subobjects (or simply denoted as $\mathcal{C}$), a normal cone with {vertex} $a \in v\mathcal{C}$ is an 'order preserving' mapping from the set $v\mathcal{C}$ (which coincides with the object set $v\mathcal{P}$ of the preorder) to the set $\mathcal{C}(c,a)$ of all morphisms in $\mathcal{C}$ from an object $c$ to codomain $a$. Normal cones can be composed via a special binary operation and all the {normal cones} in a normal category $\mathcal{C}$ form a regular semigroup $T\mathcal{C}$ known as the semigroup of normal cones in $\mathcal{C}$ . The fact that every principal ideal in a regular semigroup has an idempotent generator is captured by stipulating that every object in a normal category should have an associated idempotent normal cone.

Hence we get two normal categories, one each from the principal left and right ideals. To describe the relationship between them, a notion of {duality} is needed. Grillet defined the dual $P^*$ of a regular partially ordered set $P$ as the set of all {normal equivalence relations} on $P$. An equivalence relation was called normal, if it was induced by a normal mapping.

On the same lines, Nambooripad called a set-valued functor $H(\gamma;-)\colon\mathcal{C}\to \mathbf{Set}$ on a normal category $\mathcal{C}$ as an {$H$-functor} if it was 'induced' by a normal cone $\gamma$. He defined the \emph{normal dual} $N^*\mathcal{C}$ of a normal category $\mathcal{C}$ as the category whose object set was the set of all $H$-functors in $\mathcal{C}$. The category $N^\ast\mathcal{C}$ is a full subcategory of the category $\mathcal{C}^\ast$ where $\mathcal{C}^\ast$ is the category of all functors from $\mathcal{C}$ to $\mathbf{Set}$ \cite{mac}. Hence the morphisms in $N^*\mathcal{C}$ are natural transformations and it can be shown that the category $N^*\mathcal{C}$ forms a normal category.
Further, the Green relations in the semigroup $T\mathcal{C}$ can be completely characterized in terms of the $H$-functors and the vertices of the normal cones in $\mathcal{C}$.

Further, in a regular semigroup, the category of principal right (left) ideals is {locally isomorphic} with the normal dual of the category of principal left (right) ideals. The notion of local isomorphisms in cross-connections was initially used by Rajan \cite{loc} to study the order-isomorphisms of principal ideals of regular semigroups. He showed that local isomorphisms arise from the $\omega$-isomorphisms of the biordered set of the semigroup. This relationship between the principal ideals is abstracted as a {cross-connection} functor.

A local isomorphism between two normal categories $\mathcal{C}$ and $\mathcal{D}$ is said to be {total} if its image is a {total ideal} in $\mathcal{D}$. This is a property of a local isomorphism in regular semigroups which defines a unique local isomorphism in the opposite direction between $\mathcal{D}$ and $\mathcal{C}$.
\begin{dfn}
A {cross-connection} between two normal categories $\mathcal{C}$ and $\mathcal{D}$ is a {total local isomorphism} $\Gamma \colon \mathcal{D} \to N^\ast\mathcal{C}$ where $N^\ast\mathcal{C}$ is the {normal dual} of the category $\mathcal{C}$.
\end{dfn}
Given a cross-connection $\Gamma\colon\mathcal{D} \to N^\ast\mathcal{C}$, as discussed above, we have a unique {dual cross-connection} $\Delta \colon\mathcal{C} \to N^\ast\mathcal{D}$ such that there is a {natural isomorphism} $\chi_\Gamma$ between the bifunctors $\Gamma(-,-)$ and $\Delta(-,-)$ associated with the functors $\Gamma$ and $\Delta$. Using the natural isomorphism $\chi_\Gamma$, we can get a \emph{linking} of some normal cones $\gamma \in T\mathcal{C}$ with some normal cones $\delta \in T\mathcal{D}$.

These linked cone pairs $(\gamma,\delta)$ with respect to induced binary compositions will form a regular semigroup which is called the {cross-connection semigroup} $\tilde{S}\Gamma$ determined by $\Gamma$. If we do the construction starting with the categories of principal left and right ideals of a regular semigroup, then the cross-connection semigroup $\tilde{S}\Gamma$ is isomorphic to the regular semigroup we started with. Hence this gives a faithful representation of the semigroup as a subdirect product of $T\mathcal{C} \times (T\mathcal{D})^\text{op}$. An interested reader may refer \cite{tx,var,azeefcross,kvn,cross} for a detailed discussion on cross-connection theory.

Let $V$ be an arbitrary vectorspace over a field $K$. The singular linear transformation semigroup $Sing(V)$ is the multiplicative semigroup of all singular linear transformations on $V$. It is one of the most important regular subsemigroups of the regular monoid $\mathscr{T}_V$ of all (including non-singular) linear transformations on $V$. The full linear transformation semigroup $\mathscr{T}_V$ is one of the most important semigroups as it is a generalization of the semigroup of matrices, semigroup of operator algebras etc, and consequently has been studied extensively. 
 
The theory of cross-connections degenerates to a major extend in the case of regular monoids and hence yields very little information about the structure of the monoids. So we focus our attention on the subsemigroup $Sing(V)$, and using our approach, the analysis can be easily extended to $\mathscr{T}_V$ and even to its variants (see Section \ref{variant}). The cross-connections of $\mathscr{T}_V$ were studied in detail by Rajendran and Nambooripad \cite{raj} using a different approach where they described cross-connection semigroups arising from the structure of bilinear forms. They had restricted the discussion to finite dimensional vector spaces, but here we extend the study to arbitrary vector spaces also.
 
The structure of the article is as follows. In Section \ref{secsub}, we describe the normal categories arising from the semigroup $Sing(V)$. The category of principal left ideals of $Sing(V)$ is characterized as the Subspace category $\mathscr{S}(V)$ of proper subspaces of $V$ with linear transformations as morphisms. We show that the semigroup of all normal cones in $\mathscr{S}(V)$ is isomorphic to $Sing(V)$. We characterize the normal dual of the Subspace category in terms of the algebraic dualspace of $V$ as the Annihilator category $\mathscr{A}(V)$. When $V$ is finite dimensional, it is shown that the Annihilator category $\mathscr{A}(V)$ is isomorphic to the category $\mathscr{S}(V^*)$ of proper subspaces of $V^*$ where $V^*$ is the algebraic dual space of $V$. In Section \ref{seccxn}, we describe various cross-connections between these categories. We show how each automorphism $\theta$ on $V$ gives rise to a cross-connection between $\mathscr{A}(V)$ and $\mathscr{S}(V)$. We also prove that \emph{every} cross-connection semigroup that arises from the cross-connections is isomorphic to $Sing(V)$. In the last section, we construct some interesting semigroups arising from the cross-connections by restricting the categories. This discussion naturally leads to a useful description of the regular part of the variant of full linear transformation semigroup $\mathscr{T}_V$.

\section{Subspaces and annihilators}\label{secsub}
In this section, we discuss the normal categories arising from the semigroup $Sing(V)$. In this article, all the functions are written in the order of their composition, i.e., from left to right. A functor $F$ is said to be $v$-injective ($v$-surjective) if $F$ is injective (surjective) on the set of objects. Recall that normal categories are abstractions of principal left (right) ideals and so we will require the following well-known results regarding the Green relations in $Sing(V)$. In the sequel, we shall denote by $N_\alpha$, the null space of $\alpha$ consisting of all $y \in V$ such that $y\alpha = 0$. We denote the image of the linear transformation $\alpha$, interchangeably by $V\alpha$ and Im $\alpha$, to suit context.
\begin{lem} \label{tlx} \cite[Section 2.2]{clif}  
Let $\alpha , \beta $ be arbitrary linear transformations on $V$. 
\begin{enumerate}
\item There exists $ \: \varepsilon \in Sing(V) $ such that $\alpha\varepsilon = \beta $ if and only if $ V\alpha \supseteq V\beta $. Hence $\alpha \mathscr{L} \beta \iff V\alpha = V\beta$, i.e., Im $\alpha = $ Im $\beta$.
\item There exists $ \: \varepsilon \in Sing(V) $ such that $\varepsilon\alpha = \beta $ if and only if $ N_\alpha \subseteq N_\beta $. Hence $\alpha \mathscr{R} \beta \iff N_\alpha = N_\beta$.
\item If $\alpha \in Sing(V)$ is an idempotent, then $ V = N_\alpha \oplus V\alpha $.
\item If $V= N \oplus W$ such that $N,W \subsetneq V$, then there exists $\epsilon \in Sing(V)$ such that $N= N_\epsilon$ and $W=V\epsilon$.
\end{enumerate}
\end{lem}
It follows from the previous lemma that a description of principal left ideals of $Sing(V)$ would involve subspaces $V\alpha$ of $V$. The proper subspaces of a vectorspace $V$ with linear transformations as morphisms, in fact form a category $\mathscr{S}(V)$ which shall be called the \emph{Subspace category}. The category $\mathscr{S}(V)$ has a natural {choice of subobjects}---the one provided by subspace inclusions.

Given any linear transformation $f$ between subspaces $A$ and $B$, it has a factorisation of the form $f = quj$ where $q\colon A\to A'$ is a projection, $u=f_{|A'}$ is an isomorphism and $j=j(B',B)$ is an inclusion. Here $A'$ is a complement of the null space $N_f$ of $f$ in $A$; $q\colon A\to A'$ is the projection associated with the direct sum decomposition $N_f \oplus A' =A$, and $B' =$ Im $f = Af$. This gives the normal factorization of $f$ and the linear transformation $qu\colon A \to \text{Im } f$ is called the \emph{epimorphic component} $f^\circ$ of $f$.
 
Let $v\mathscr{S}(V)$ and $\mathscr{S}(V)$ denote the object set and the morphism set of the category $\mathscr{S}(V)$ respectively. Given any subspace $D \subseteq V$, we associate a mapping $\sigma\colon v\mathscr{S}(V) \to \mathscr{S}(V)$ with the following properties.
\begin{enumerate}
\item  For each subspace $A$ of $V$, $\sigma(A)$ is a linear transformation from $  A$ to $ D$ and whenever $A \subseteq B$, $j(A,B)\sigma(B) =\sigma(A)$.
\item For some subspace $C$ of $V$, $\sigma(C)\colon C \to D$ is an isomorphism.
\end{enumerate}
This collection of morphisms $\{\sigma(A) \: : \: A\in v\mathscr{S}(V)\}$ is called a \emph{normal cone} $\sigma$ with vertex $D$ in the category $\mathscr{S}(V)$. We stress that a normal cone is \emph{not} defined by its vertex; there may exist many normal cones sharing the same vertex.

Let $\gamma$ be a normal cone with vertex $C$ and $f\colon  C\to D$ an onto linear transformation in $\mathscr{S}(V)$. Then we can see that $\gamma*f\colon  A\mapsto \gamma(A) f$ gives a normal cone in $\mathscr{S}(V)$ with vertex $D$. Let $\gamma$, $\delta$ are two normal cones in $\mathscr{S}(V)$ with vertices $C$ and $D$ respectively, then $(\delta(C))^\circ$, the epimorphic component of the morphism $\delta(C)$ is an onto linear transformation from $C$ to Im $\delta(C)$. So, we may compose the normal cones $\gamma$ and $\delta$ as follows. 
\begin{equation} \label{eqnsg}
\gamma \cdot \delta := \gamma \ast (\delta(C))^\circ
\end{equation} 
Then it can be seen that $\gamma \cdot \delta$ is a normal cone with vertex  Im $\delta(C) \: \subseteq D$. Further, the normal cone $\sigma$ will be called an {idempotent} normal cone, if $\sigma(D) =1_D$.

Let $u\colon  V \to D $ be a linear transformation such that $(x)u = x \quad \forall x \in D$. For any subspace $A \subseteq V$, define 
$$\sigma(A) := u_{|A} \colon  A \to D.$$
Then $\sigma$ is an idempotent normal cone with $\sigma(D) = 1_D$ and hence $\mathscr{S}(V)$ is a normal category. The set of all normal cones in $\mathscr{S}(V)$ under the binary operation defined in (\ref{eqnsg}) forms a regular semigroup $T\mathscr{S}(V)$, and this is the {semigroup of normal cones} in $\mathscr{S}(V)$.

Given an element $\alpha$ in $Sing(V)$, we have a normal cone $\rho^\alpha$ in $\mathscr{S}(V)$ called a {principal cone} determined by $\alpha$. Observe that by Lemma \ref{tlx}, for each idempotent $e \in Sing(V)$, the principal left ideal $ (Sing(V))e$ is identified with the subspace Im $e$. So the principal cone $\rho^\alpha$ is defined, for each subspace Im $e$, as follows. 
$$ \rho^\alpha(\text{Im } e) := e\alpha_{|\text{Im e}}.$$
The morphism $e\alpha_{|\text{Im e}}$ maps an arbitrary element $\theta$ in the principal left ideal $ (Sing(V))e$ to $\theta e\alpha$ in the principal left ideal $ (Sing(V))\alpha$. This collection of morphisms clearly forms a normal cone in $\mathscr{S}(V)$ with vertex Im $\alpha$. It is known that the map $\alpha \mapsto\rho^\alpha$ is a homomorphism from a given regular semigroup to its semigroup of normal cones. Although it is an isomorphism in the case of regular monoids (which in fact leads to a degeneration of the analysis), this map need not be injective or onto, in general. For instance, the map fails to be injective in the case of rectangular bands and it is not onto in arbitrary completely simple semigroups \cite{amth,css}. Section \ref{variant} gives another instance where the map is not onto. But in the case of $Sing(V)$, every normal cone $\sigma$ with vertex $D$ in $\mathscr{S}(V)$ defines a linear transformation $\alpha\colon  X \to D$ such that $\sigma =\rho^\alpha$ as follows. If $B$ is a basis of $V$, let 
\begin{equation}\label{cone}
(b) \alpha := (b)\sigma(\langle b \rangle) \text{ for all } b \in B
\end{equation} 
where $\sigma(\langle b \rangle)$ is the component of $\sigma$ at the subspace $\langle b \rangle \in v \mathscr{S}(V)$. Since a linear transformation on $V$ is completely determined by the values on its basis elements, $\alpha$ is a linear transformation on $V$ such that $\sigma =\rho^\alpha$.

Thus every normal cone in $\mathscr{S}(V)$ is a principal cone and hence the map is onto. Further suppose $\rho^\alpha = \rho^\beta$. Then for any $b \in B$, since $\rho^\alpha (\langle b \rangle) = \rho^\beta (\langle b \rangle) $, we see that $ b\alpha = b \beta$ for all $b \in B$, and hence $\alpha = \beta$. Summarising, we have the following theorem.
\begin{thm}\label{thm1}
$\mathscr{S}(V)$ is a normal category and $T\mathscr{S}(V)$ is isomorphic to $Sing(V)$.
\end{thm}  
Now we proceed to characterize the {normal dual} $N^*\mathscr{S}(V)$ of the Subspace category $\mathscr{S}(V)$. The objects of the {normal dual} $N^*\mathscr{S}(V)$ are $H$-functors $H(\gamma;-)\colon \mathscr{S}(V) \to \mathbf{Set}$ where $\gamma$ is an idempotent normal cone in $\mathscr{S}(V)$. By Theorem \ref{thm1}, we have an identification of the normal cones in $T\mathscr{S}(V)$ with linear transformations in $Sing(V)$. Hence, we define an $H$-functor in $Sing(V)$ as follows. Let $e$ be an idempotent transformation in $Sing(V)$. For each $A \in v\mathscr{S}(V)$ and for each morphism $g\colon A\to B$ in $\mathscr{S}(V)$, $ H(e;-) \colon  \mathscr{S}(V) \to \mathbf{Set}$ is a functor such that 
\begin{subequations} 
\begin{align}
H(e;A)= &\{eh : \: h:\text{Im }e\to A \} \text{ and } \label{eqnH1}\\
H(e;g) \colon H(e;A) &\to H(e;B) \text{ given by }eh \mapsto ehg.\label{eqnH2}
\end{align}
\end{subequations}
The following lemma characterizing $H$-functors is pretty straight forward but is crucial in the sequel.
\begin{lem}\label{leml2}
Let $e \in  Sing(V) $ and $A\subseteq V$. Then
$$ H (e ; {A}) = \{  a \in Sing(V) \: : \: N_a \supseteq N_e \text{ and }Va\subseteq A  \}.$$
If $g \colon  {A}\to{B} $ then $ H (e ; g ) \colon  H(e;{A})  \to  H(e;{B})  $ is given by $ a \mapsto  ag $ for $a \in H(e;A)$. 
\end{lem}
The above lemma tells us that the $H$-functors $H(e;-)$ in $\mathscr{S}(V)$ is completely determined by the null space $N_e$ of $e$. This inspires us to describe the normal dual $N^*\mathscr{S}(V)$ using the algebraic duality of $V$. Also recall that the principal right ideals of $Sing(V)$ is characterized by the null spaces of $V$.

Given an arbitrary vectorspace $V$ over the field $K$, a {linear functional}\index{linear functional} on $V$ is a linear transformation from $V$ to $K$. The set of all the linear functionals on the vectorspace $V$ forms a vectorspace $V^*$ and is known as the {algebraic dual space} of $V$. Given a linear transformation $u \colon  A \to B$, there is a unique linear transformation $u^* \colon  B^* \to A^*$ called the {transpose}\index{transpose} of $f$ given by 
\begin{equation}\label{eqntrans}
(\alpha)u^* = u\alpha \text{ for all } \alpha \in B^*.
\end{equation}  
Now we define a new category $\mathscr{A}(V)$ called the \emph{Annihilator category} whose objects are $A^\circ$ where $A^\circ \:=\: \{ f \in V^* : vf=0 \text{ for all } v \in A \}$ is the annihilator\index{annihilator} of $A$, for a proper subspace $A$ of $V$.

To define morphisms in $\mathscr{A}(V)$, we need to analyse the morphisms in $N^*\mathscr{S}(V)$, which are in fact natural transformations between $H$-functors. It has been shown in \cite[III, Lemma 6]{cross} that an $H$-functor $H(e;-)$ is a \emph{representable functor} with Im $e$, as the representing object. 

This means that there is a natural isomorphism $\eta_e\colon H(e;-) \to \mathcal{S}(V)(\text{Im } e,-)$, where $\mathcal{S}(V)(\text{Im } e,-)$ is the covariant hom-functor in the category $\mathcal{S}(V)$ determined by the object Im $e$. For a subspace $A\in v\mathcal{S}(V)$ and for $h\colon\text{Im }e \to A$ in $\mathcal{S}(V)$, the component $\eta_e(A)$ of the natural isomorphism $\eta_e$, is given by $\eta_e(A): eh\mapsto h$. Observe here that $eh \in H(e;A)$ and $h$ is a linear transformation in $\mathcal{S}(V)(\text{Im } e,A)$. Compare with (\ref{eqnH1}) and Lemma \ref{leml2} above.

Further, by \emph{Yoneda lemma}, natural transformations between the hom-functors $\mathcal{S}(V)(\text{Im } e,-)$ and $\mathcal{S}(V)(\text{Im } f,-)$ are in one-to-one correspondence with the morphisms (in the reverse direction) between the representing objects, Im $f$ and Im $e$. So, as in the proof of \cite[III, Lemma 21]{cross}, we can characterize the set of morphisms between $ H(e;-) $ and $ H({f} ;-)$ in $N^*\mathscr{S}(V)$ using the set of morphisms between Im $f$ and Im $e$ in the category $\mathscr{S}(V)$. Hence we have the following lemma.
\begin{lem}
Let $\sigma$ be a natural transformation between $ H(e;-) $ and $ H({f} ;-)$, then there is a unique $u\in f(Sing(V))e$ such that $u_{|\text{Im }f}\colon  \text{Im }f \to \text{Im }e$ and the component of $\sigma$ at a subspace $A= \text{Im }h$ is given by $\sigma(A)\colon eh\mapsto fuh$.
\end{lem}
The above lemma may also be illustrated elegantly using the following commutative diagram.

\begin{equation*}\label{Homf}
\xymatrixcolsep{2pc}\xymatrixrowsep{3pc}\xymatrix
{
 H(e;-) \ar[rr]^{\eta_e} \ar[d]_{\sigma}  
 && \mathscr{S}(V)(\text{Im }e,-) \ar[d]^{\mathcal{L}(u_{|\text{Im }f},-)} & \text{Im }e\\       
 H(f;-) \ar[rr]^{\eta_{f}} && \mathscr{S}(V)(\text{Im }f,-) & \text{Im }f \ar[u]_{u_{|\text{Im }f}}
}
\end{equation*}

The isomorphism $\eta_e(A)$ maps an element $eh \in H(e;A)$ to $h$ in $ \mathscr{S}(V)(\text{Im }e,A)$. Then $\mathcal{L}(u_{|\text{Im }f},A)$ maps it to $u_{|\text{Im }f}h$. Observe $u_{|\text{Im }f}h \in \mathscr{S}(V)(\text{Im }f ,A) $. Since $\eta_{f}(A)$ is an isomorphism, it is invertible. So $\eta_{f}(A)^{-1}$ maps $u_{|\text{Im }f}h \mapsto fuh$. Thus by diagram chasing, $\sigma(A)=\eta_e(A)\mathcal{L}(u_{|\text{Im }f},A)\eta_{f}(A)^{-1}:uh \mapsto fuh$ and hence the lemma.

In the context of dual spaces, we can simplify the above lemma further to get a nice description of morphims in $\mathscr{A}(V)$ as follows. First, since $N_e \subseteq V$, $(N_e)^\circ \in \mathscr{A}(V)$. If $\alpha \in (N_e)^\circ$, then $(N_e)\alpha = 0$. As $u \in f(Sing(V))e$ and using Lemma \ref{tlx}, we see that $N_f \subseteq N_u$, and $xu\alpha = 0$ for all $x \in N_f$. That is $(N_f)u\alpha = 0$ and $u\alpha \in (N_f)^\circ$. So $u^* \colon  \alpha \mapsto u\alpha$ is a well-defined linear transformation from $(N_e)^\circ$ to $(N_f)^\circ$. Hence, for $u\in f(Sing(V))e$, a morphism $u^*$ in $\mathscr{A}(V)$ from $(N_e)^\circ$ to $(N_f)^\circ$ is defined as follows. 
$$u^*\colon \alpha \mapsto u\alpha\text{ for all }\alpha\in (N_e)^\circ.$$
Since $V^*$ is a vectorspace, $\mathscr{S}(V^*)$ is a normal category. We can see that $\mathscr{A}(V)$ is a subcategory of $\mathscr{S}(V^*)$. Now we can define a functor $P \colon  N^\ast \mathscr{S}(V) \to \mathscr{S}(V^*)$  as follows.
\begin{equation} \label{eqnPl}
P (H(e;-)) = (N_e)^\circ \quad\text{and}  \quad
P(\sigma) = u^*
\end{equation}
where $(N_e)^\circ$ is the annihilator of the null space of $e$ and $\sigma (Vh)\colon eh\mapsto fuh$ where $u \in f(Sing(V))e$. We can see that $P$ is well-defined, inclusion preserving, $v$-injective and faithful. Since Im $P = \mathscr{A}(V)$, we have the following results.
\begin{pro}\label{produal}
The normal dual $N^*\mathscr{S}(V)$ of the Subspace category is isomorphic to $\mathscr{A}(V)$. Dually, $N^*\mathscr{A}(V)$ is isomorphic to $\mathscr{S}(V)$. The semigroup $T\mathscr{A}(V)$ of normal cones in $\mathscr{A}(V)$ is isomorphic to $Sing(V)^\text{op}$.
\end{pro}
Although in general, the functor $P$ need not be $v$-surjective or full, it is both $v$-surjective and full when $V$ is finite dimensional. In that case, we see that the normal dual of the Subspace category of the vectorspace $V$ is isomorphic to the Subspace category of the algebraic dual space $V^*$ of $V$ as described in the following theorem. Here the cross-connection duality coincides with the algebraic duality of vectorspaces.
\begin{thm} \label{thmplx}
Let $V$ be a finite dimensional vectorspace over $K$, then $N^\ast \mathscr{S}(V)$ is isomorphic as a normal category to $\mathscr{S}(V^*)$.
\end{thm}

\section{Cross-connections}\label{seccxn}
Summarising, we have seen that the two categories involved in the construction of $Sing(V)$ are the Subspace category $\mathscr{S}(V)$ and the Annihilator category $\mathscr{A}(V)$, and they are mutual normal duals. Now we proceed to describe \emph{all} the cross-connections between them, and construct the semigroups arising from these categories. We will need the following definitions from the general cross-connection theory in the sequel.
\begin{dfn}\label{ideal}
Let $\mathcal{C}$ be a normal category. Then an \emph{ideal}\index{ideal!of a category} $\langle c \rangle$ of $\mathcal{C}$ is the full subcategory of $\mathcal{C}$ whose objects are subobjects of $c$ in $\mathcal{C}$. It is called the principal ideal generated by $c$.
\end{dfn}
\begin{dfn}\label{lociso}
Let $\mathcal{C}$ and $\mathcal{D}$ be normal categories. Then a functor $F\colon  \mathcal{C} \to \mathcal{D}$ is said to be a \emph{local isomorphism}\index{local isomorphism} if $F$ is inclusion preserving, fully faithful and for each $c \in v\mathcal{C}$, $F_{|\langle c \rangle}$ is an isomorphism of the ideal $\langle c \rangle$ onto $\langle F(c) \rangle$.
\end{dfn}
\begin{dfn}\label{cxn}
A \emph{cross-connection}\index{cross-connection} from $\mathcal{D}$ to $\mathcal{C}$ is a triplet $(\mathcal{D},\mathcal{C};{\Gamma})$ where $\Gamma\colon  \mathcal{D} \to N^\ast\mathcal{C}$ is a local isomorphism such that for every $c \in v\mathcal{C}$, there is some $d \in v\mathcal{D}$ such that $c \in M\Gamma(d)$.
\end{dfn}
The M-set\index{M-set} associated with a cone $\sigma$ in $\mathscr{S}(V)$ (also written $MH(\sigma;-)$) is given by
$$ M\sigma = \{ A \in v\mathscr{S}(V) \: |\: \sigma(A)\text{ is an isomorphism} \}.$$
In the case of the Subspace category, the M-set may be characterized as follows.
\begin{pro}\label{proml}
The M-set\index{M-set} of the cone $\rho^e$ is given by $M{\rho^e} = MH(e;-) = M((N_e)^\circ)= \{ A\subseteq V : A \oplus N_e = V \} $.
\end{pro}
From the previous discussion, $(\mathscr{A}(V),\mathscr{S}(V),\Gamma)$ is a cross-connection if $\Gamma\colon  \mathscr{A}(V) \to \mathscr{A}(V)$ is a local isomorphism such that for every $A \in v\mathscr{S}(V)$, there is some $Y \in v\mathscr{A}(V)$ such that $A \in M(\Gamma(Y))$. Given a cross-connection $\Gamma \colon  \mathscr{A}(V) \to \mathscr{A}(V)$ with a dual $\Delta \colon  \mathscr{S}(V) \to \mathscr{S}(V)$, we have two associated bifunctors\index{functor!bifunctor} $\Gamma \colon  \mathscr{S}(V) \times \mathscr{A}(V) \to \bf{Set}$ and $\Delta \colon  \mathscr{S}(V) \times \mathscr{A}(V) \to \bf{Set}$ such that for all $(A,Y) \in v\mathscr{S}(V) \times v\mathscr{A}(V)$ and $(f,w^*)\colon (A,Y) \to (B,Z)$
\begin{subequations}\label{eqnglx}
\begin{align}
\Gamma (A,Y)\:&=\: \{ \alpha \in Sing(V) : V \alpha \subseteq A \text{ and } (N_\alpha)^\circ \subseteq \Gamma(Y) \}\\
\Gamma (f,w^*)\:& \colon \: \alpha \mapsto (y\alpha) f \:=\: y(\alpha f)
\end{align}
\end{subequations}
where $y$ is given by $y^* = \Gamma(w^*)$, and 
\begin{subequations}\label{eqndlx}
\begin{align}
\Delta (A,Y)\:&=\: \{ \alpha \in Sing(V) : V\alpha \subseteq \Delta(A) \text{ and } (N_\alpha)^\circ \subseteq Y \}\\
\Delta (f,w^*)\:&\colon \: \alpha \mapsto (w\alpha) g \:=\: w(\alpha g)
\end{align}
\end{subequations} 
where $g = \Delta(f)$.

From the above discussion, it may be observed that every isomorphism between the categories will indeed be a cross-connection. In order to describe isomorphisms between the subspaces of a vectorspace, naturally our attention shifts to the automorphisms of $V$.

Recall that an {automorphism}\index{automorphism} $\theta$ of a vectorspace $V$ is an isomorphism from $V$ onto itself. In the sequel, since $\theta$ often represents a functor, we write the argument on the right side of the automorphism $\theta$ also, for uniformity, but the order of composition is from left to right as earlier. 

First observe that the automorphism $\theta$ will map a proper subspace $A$ of $V$ isomorphically to its image. So define a functor $\Delta_{\theta} \colon  \mathscr{S}(V) \to \mathscr{S}(V)$ as follows. For a proper subspace $A$ of $V$ and for $f: A\to B$ in $\mathscr{S}(V)$,
\begin{equation}\label{eqndcrossrpl}
\Delta_{\theta} (A) \:=\: \theta(A) \text{ and } \Delta_{\theta}(f) = \theta(f) = \theta^{-1} f \theta.
\end{equation} 
Then the functor $\Delta_\theta$ is a normal category isomorphism on the Subspace category $\mathscr{S}(V)$ and so $(\mathscr{S}(V), \mathscr{A}(V) ; \Delta_{\theta})$ is a cross-connection.
 
Since $\theta$ is an automorphism on $V$, the transpose $\theta^\ast : V^* \to V^*$ is also an automorphism. Hence define a functor $\Gamma_\theta$ on $\mathscr{A}(V)$ as follows. For a proper subspace $Y$ of $V^*$, and for $w^*\colon Y\to Z$ in $\mathscr{A}(V)$,
\begin{equation}\label{eqncrossrpl}
\Gamma_{\theta} (Y) \:=\: \theta^*(Y) \text{ and } \Gamma_{\theta}(w^*) = (\theta^{-1})^* w^* (\theta^*) = (\theta w \theta^{-1})^*.
\end{equation} 
Then we can see that the functor $\Gamma_\theta$ is a normal category isomorphism on $\mathscr{A}(V)$ and so $(\mathscr{A}(V),\mathscr{S}(V); \Gamma_{\theta})$ is also a cross-connection. 

From the above cross-connections, we obtain two bifunctors $\Gamma_{\theta}(-,-)$ and $\Delta_{\theta } (-,-)$ respectively as in (\ref{eqnglx}) and (\ref{eqndlx}). Now if we define a natural transformation $\chi_{\Gamma_\theta}\colon \Gamma_{\theta }(-,-) \to \Delta_{\theta }(-,-)$ as $\chi_{\Gamma_\theta}(A,Y)\colon  \alpha \mapsto \theta^{-1}\alpha\theta$, then for $(f,w^*) \in \mathscr{S}(V) \times \mathscr{S}(V^*)$ with $f \colon  A \to B$ and $w^* \colon  Y \to Z$, the following diagram commutes.
\begin{equation*}\label{duality}
\xymatrixcolsep{4pc}\xymatrixrowsep{3pc}\xymatrix
{
 \Gamma_{\theta }(A,Y) \ar[r]^{\chi_{\Gamma_\theta}(A,Y)} \ar[d]_{\Gamma_{\theta }(f,w^\ast)}  
 & \Delta_{\theta }(A,Y) \ar[d]^{\Delta_{\theta }(f,w^\ast)} \\       
 \Gamma_{\theta }(B,Z) \ar[r]^{\chi_{\Gamma_\theta}(B,Z)} & \Delta_{\theta }(B,Z) 
}
\end{equation*}
Also since the map $\chi_{\Gamma_\theta}(A,Y)$ is a bijection of the set $\Gamma_{\theta }(A,Y)$ onto the set $\Delta_{\theta }(A,Y)$, $\chi_{\Gamma_\theta}$ is a natural isomorphism. Hence the cross-connection $\Delta_\theta$ is the dual of the cross-connection $\Gamma_\theta$ and $\chi_{\Gamma_\theta}$ is the duality associated with $\Gamma_\theta$.
 
Consequently we see that $\alpha \in T\mathscr{S}(V)$ is linked to $\beta\in T\mathscr{A}(V)$ via the duality $\chi_{\Gamma_\theta}$ if and only if $ \beta = \theta^{-1}\alpha\theta$. So the cross-connection semigroup is given by
$$ \tilde{S}\Gamma_\theta = \:\{\: (\alpha,\theta^{-1}\alpha\theta) \: \text{ such that } \alpha \in Sing(V)\}.$$
Hence clearly $\tilde{S}\Gamma_\theta$ is isomorphic to $Sing(V)$.

Now given an arbitrary cross-connection $\Gamma$ from $\mathscr{A}(V)$ to $\mathscr{S}(V)$ with the dual cross-connection $\Delta$ from $\mathscr{S}(V)$ to $\mathscr{A}(V)$, define a mapping $\theta:V \to V$ as:
$$\theta (b)= x \text{ such that } \Delta( \langle b \rangle) = \langle x \rangle$$
for all $b \in B$, where $B$ is a basis of $V$. Then $\theta$ will clearly be a linear transformation on $V$. Using the cross-connection properties of $\Delta$, the linear transformation $\theta$ can be shown to be an automorphism on $V$. Arguing on similar lines as in the case of normal cones in $\mathscr{S}(V)$, we see that $\Gamma = \Gamma_\theta$ and $\Delta = \Delta_{\theta}$. Since the cross-connection semigroup $ \tilde{S}\Gamma_\theta$ is isomorphic to $Sing(V)$, we have the following theorem.
\begin{thm}
Every cross-connection semigroup arising from the cross-connections between the categories $\mathscr{A}(V)$ and $\mathscr{S}(V)$ is isomorphic to $Sing(V)$.
\end{thm}

\section{Variants}\label{variant}
So far we have described the normal categories $\mathscr{S}(V)$ and $\mathscr{A}(V)$, and studied how they can be \emph{connected} to construct semigroups. The last theorem shows that they do not give rise to any 'new' semigroup. Observe that by adjoining $V$ to $\mathscr{S}(V)$, and $V^*$ to $\mathscr{A}(V)$, we can characterize the categories in the {full} linear transformation semigroup $\mathscr{T}_V$. And using similar arguments (in fact easier, since the \emph{largest object}---$V$ should map to $V$ under cross-connections), we can show that every semigroup that arises is isomorphic to $\mathscr{T}_V$. In this section, we aim to construct the regular part of the {non-trivial} variants of $\mathscr{T}_V$ by {cross-connecting} certain subcategories of $\mathscr{S}(V)$ and $\mathscr{A}(V)$. In the sequel, we will be denoting the categories arising from $\mathscr{T}_V$ also using $\mathscr{S}(V)$ and $\mathscr{A}(V)$, and whether the largest object is adjoined or not will be clear from the context.

For an arbitary linear transformation $\theta$, let $\mathscr{T}_V^\theta = (\mathscr{T}_V,\ast)$ be the variant of the full linear transformation semigroup with the binary composition $\ast$ defined as follows. 
$$ \alpha \ast \beta = \alpha \cdot \theta \cdot \beta \quad \text{ for } \alpha, \beta \in \mathscr{T}_V.$$
Variant of arbitrary semigroups was initially studied by Hickey \cite{hickey}, and later by Khan and Lawson \cite{khan}. Khan and Lawson, among other things, had showed that the regular elements of the variant of a regular semigroup form a regular subsemigroup. The variants of the full transformation semigroup have been of recent interest following the works of Tsyaputa \cite{tsya}, Ganyushkin and Mazorchuk \cite{cfts}, Kemprasit \cite{kempra}, Dolinka and East \cite{igd} etc.

Dolinka and East \cite{igd,igd2} have studied the variants of finite full transformation semigroup and semigroup of rectangular matrices in great detail. They have described the structure of the regular part Reg$(\mathscr{T}_X^\theta)$ using the map $\alpha\mapsto (\alpha\theta,\theta\alpha)$ starting from the right and left translations. We will also give a similar description for Reg$(\mathscr{T}_V^\theta)$ as $\alpha \mapsto (\theta\alpha,\alpha\theta)$, but using cross-connections. Observe the change in the order of components although in both descriptions, the first coordinates arose from right translations. Also our description comes from a much more general theory, and so we expect that it will shed some light on the problem of extending the structure theorem to a more general class of variants.

First, observe that the cross-connection semigroup $ \tilde{S}\Gamma_\theta$ described in the last section (by extending the discussion to the full linear transformation semigroup $\mathscr{T}_V$) refers to the cross-connection arising from the variant semigroup $\mathscr{T}_V^\theta$ where $\theta$ is an automorphism. 
If we define $\phi\colon  \mathscr{T}_V^\theta \to \tilde{S}\Gamma_\theta$ as $\alpha \mapsto (\theta\alpha,\alpha\theta)$, then 
$$(\alpha \ast \alpha')\phi = (\alpha\theta\alpha')\phi = (\theta\alpha\theta\alpha',\alpha\theta\alpha'\theta ) = (\theta\alpha, \alpha\theta)\circ(\theta \alpha',\alpha'\theta) = (\alpha)\phi \circ (\alpha')\phi.$$ 
Further it can be seen that $\phi$ is an isomorphism. But $\phi$ is effectively the map $\alpha \mapsto (\rho^{\theta\alpha},\lambda^{\alpha\theta})$ where $\rho^{\theta\alpha}$ and $\lambda^{\alpha\theta}$ are precisely the principal cones in $T\mathscr{S}(V)$ and $T\mathscr{A}(V)$ respectively, which arise from the semigroup $\mathscr{T}_V^\theta$. Hence the same construction may be extended to the variant semigroup for an arbitrary $\theta$. But when $\theta$ is not an automorphism, the map is injective only in the regular part and so we can only recover Reg$(\mathscr{T}_V^\theta)$. 

Now we define categories $\mathscr{R}(V)$ and $\mathscr{B}(V)$ as normal full-subcategories of $\mathscr{S}(V)$ and $\mathscr{A}(V)$ respectively, as follows:
$$v\mathscr{R}(V) = \{ A : \: A \subseteq (N_\theta)^c \} \text{ and } v\mathscr{B}(V) = \{ A^\circ : \: N_\theta \subseteq A \}.$$
Here $(N_\theta)^c$ is a complementary subspace to the null space of $\theta$. Then the semigroup of normal cones in $\mathscr{R}(V)$ is given by $T\mathscr{R}(V) = \{ f\in \mathscr{T}_V : \text{Im }f \subseteq (N_\theta)^c\} $, and the semigroup of normal cones in $\mathscr{B}(V)$ is $T\mathscr{B}(V) = \{ f\in \mathscr{T}_V : N_\theta \subseteq N_f\} $. Clearly $T\mathscr{R}(V)$ and $T\mathscr{B}(V)$ are regular subsemigroups of $\mathscr{T}_V$ and $\mathscr{T}_V^\text{op}$ respectively.

If $\mathscr{R}(V)$ is \emph{big} enough (and that depends on $\theta$), the normal dual $N^*\mathscr{R}(V)$ will be $\mathscr{A}(V)$, else a proper normal full-subcategory of it. Similarly the normal dual $N^*\mathscr{B}(V)$ will be a normal full-subcategory of $\mathscr{S}(V)$.

Now we describe how these categories may be {connected}. For $A \in v\mathscr{R}(V)$, and for $f\colon A\to B$ in $\mathscr{R}(V)$, 
$$\Delta(A) = \theta(A) \text{ and } \Delta(f) = (\theta_{|A})^{-1} f \theta.$$ 
Since $A \subseteq (N_\theta)^c$, the functor $\Delta$ is a {local isomorphism} from $\mathscr{R}(V)$ to $N^*\mathscr{B}(V)$, and hence a cross-connection. Similarly, for $\alpha \in A^\circ \in v\mathscr{B}(V)$, $N_\theta\alpha = 0$. So $\theta^*\colon \alpha \mapsto \theta . \alpha$ is a local isomorphism from $A^\circ$ to $\theta^*(A^\circ)$. So for $A^\circ \in v\mathscr{B}(V)$ and for $w^{*} \colon  A^\circ \to B^\circ$, define $\Gamma$ as follows: 
$$\Gamma (A^\circ) \:=\: \theta^*(A^\circ) \text{ and } \Gamma(w^*) = (\theta^*_{|A^\circ})^{-1} w^* \theta^*.$$
Then the transpose $\theta^*$ induces a dual cross-connection $\Gamma$ from the category $\mathscr{B}(V)$ to $N^*\mathscr{R}(V)$.

Imitating the construction as in the previous section (with necessary changes), we obtain a regular semigroup. We can show that the cross-connection semigroup that arises from the above cross-connection is isomorphic to the semigroup $Reg(\mathscr{T}_V^\theta)$ of all regular elements in $\mathscr{T}_V^\theta$. Thus we have a representation of $Reg(\mathscr{T}_V^\theta)$ as a subdirect product of $\mathscr{T}_V \times \mathscr{T}_V^\text{op}$ given by $\alpha \mapsto (\theta\alpha,\alpha\theta)$. Summarising, we have the following.
\begin{thm}
The cross-connection semigroup $(\mathscr{B}(V),\mathscr{R}(V); \Gamma)$ is isomorphic to the semigroup $Reg(\mathscr{T}_V^\theta)$ of all regular elements of $\mathscr{T}_V^\theta$.
\end{thm}
The above discussion shows that the principal ideals in $Reg(\mathscr{T}_V^\theta)$ are the categories $\mathscr{B}(V)$ and $\mathscr{R}(V)$, and they are {cross-connected} with each other via the sandwich element $\theta$. Observe here that only $\{ \theta\alpha : \alpha \in\text{ Reg}(\mathscr{T}_V^\theta) \}$ are the principal cones in $T\mathscr{R}(V)$, and so there are non-principal normal cones in $\mathscr{R}(V)$; similarly in the category $\mathscr{B}(V)$. Also observe that the functors $\Gamma$ and $\Delta$ are {proper local isomorphisms} which are not category isomorphisms.
 
In $Reg(\mathscr{T}_V^\theta)$, all the subtleties of cross-connection theory play out nicely; reiterating the necessity of such a sophisticated theory in describing the ideal structure of semigroups. This in turn suggests that whenever a semigroup with a complicated ideal structure arises, it is indeed worth taking the risk of 'crossing' into cross-connections! 

How we may extend this approach to the entire variant semigroup is an open problem, and the author believes that a solution to this problem may shed some light into the much more general problem of the cross-connection construction of arbitrary semigroups. It is worth mentioning here that Nambooripad himself has made some advances in this direction using certain categories called \emph{set-based categories} \cite{newcross,newcross1}.

\section*{Acknowledgement}
The author is grateful to A. R. Rajan, University of Kerala, Thiruvananthapuram, India for several fruitful discussions during the preparation of this article.\\
The author is grateful to the editor M. V. Volkov, Ural Federal University, Russia for his keen interest and suggestions which helped improve the manuscript considerably.  

\bibliographystyle{plain}
\bibliography{SF0}

\end{document}